\documentclass[oneside,english]{amsart}
\usepackage[T1]{fontenc}
\usepackage[latin9]{inputenc}
\usepackage{babel}
\usepackage{array}
\usepackage{float}
\usepackage{multirow}
\usepackage{amstext}
\usepackage{amsthm}
\usepackage{amssymb}
\usepackage{graphicx}
\usepackage[unicode=true,pdfusetitle,
 bookmarks=true,bookmarksnumbered=false,bookmarksopen=false,
 breaklinks=false,pdfborder={0 0 1},backref=false,colorlinks=false]
 {hyperref}
\usepackage{breakurl}

\makeatletter

\providecommand{\tabularnewline}{\\}
\floatstyle{ruled}
\newfloat{algorithm}{tbp}{loa}
\providecommand{\algorithmname}{Algorithm}
\floatname{algorithm}{\protect\algorithmname}

\numberwithin{equation}{section}
\numberwithin{figure}{section}
\theoremstyle{plain}
\newtheorem{thm}{\protect\theoremname}
  \theoremstyle{plain}
  \newtheorem{cor}[thm]{\protect\corollaryname}

\usepackage{algorithm,algpseudocode}
\usepackage{pgf,tikz}
\usetikzlibrary{arrows}

\@ifundefined{showcaptionsetup}{}{%
 \PassOptionsToPackage{caption=false}{subfig}}
\usepackage{subfig}
\makeatother

  \providecommand{\corollaryname}{Corollary}
\providecommand{\theoremname}{Theorem}

\begin{document}

\title{Improvements on Spectral Bisection}

\author{Israel Rocha }

\thanks{rocha@cs.cas.cz. Institute of Computer Science, Czech Academy of
Sciences. Pod Vodarenskou vezi 271/2. 18207 Praha 8, Czech Republic. }
\begin{abstract}
We investigate combinatorial properties of certain configurations
of a graph partition which are related to the minimality of a cut.
We show that such configurations are related to the third eigenvector
of the Laplacian matrix. It is well known that the second eigenvector
encodes structural information, and that can be used to approximate
a minimum bisection. In this paper, we show that the third eigenvector
carries structural information as well. We then provide a new spectral
bisection algorithm using both eigenvectors. The new algorithm is
guaranteed to return a cut that is smaller or equal to the one returned
by the classic spectral bisection. Also, we provide a spectral algorithm
that can refine a given partition and produce a smaller cut.
\end{abstract}


\maketitle

\section{Introduction}

The classic problem of finding a minimum cut of a graph is known to
be NP-hard. Nevertheless, the problem has direct applications in VLSI
design, data-mining, finite elements and communication in parallel
computing, etc. In practice, given the importance of the problem,
the solution is generally approximated using heuristic algorithms.
The problem is to separate the vertices of a graph in two parts, such
that the number of edges connecting vertices in different parts is
minimized. Such partition, also known as a cut, is called a balanced
cut or a bisection whenever both parts have the same size.

In many applications it is desired to obtain the smallest possible
cut at a cost of having a partition that is not balanced, but acceptable
in the sense both parts have almost the same size. However, even for
those cases efficient algorithms that approximates balanced cuts up
to a constant factor do not exist. In fact, this approximation problem
is NP-hard \cite{BJ}.

Spectral techniques are well-known approaches to this problem and
they have its roots in the work of Fiedler \cite{Fiedler} and Donath
and Hoffman \cite{DonathHoffman,DonathHoffman2}. These spectral methods
are known to provide good answers, and they are broadly used in several
problems \cite{Hendrickson2,Pothen,Pothen2}. Spectral partitioning
algorithms recover global structural information and connectivity
of a graph by means of the eigenvector of the second eigenvalue of
the Laplacian matrix of the graph.

In \cite{SpiemanTen}, Spielman and Ten provided a recursive spectral
bisection algorithm and showed that spectral partitioning methods
work well on bounded-degree planar graphs. In \cite{Guattery} Guattery
and Miller, perform an analysis of the quality of the separators produced
by such methods. Papers \cite{SpiemanTen} and \cite{Guattery} discuss
the difference between guarantees on the size of a balanced cut versus
its optimality. Hendrickson and Leland \cite{Hendrickson2} extend
the spectral approach to partition a graph into four or eight parts
by using multiple eigenvectors. 

In this paper, instead of using structural information provided by
multiple eigenvectors to partition graphs into multiple parts, we
develop an approach that uses multiple eigenvectors to create a bisection
of the graphs. It is well known that the second eigenvector encodes
structural information, and that can be used to approximate a minimum
bisection. In this paper, we show that the third eigenvector carries
structural information as well, which enable us to apply that information
in the bisection problem. We then provide a new spectral bisection
algorithm using both eigenvectors.

From a more general perspective, there are several heuristics for
the graph partitioning problem, and they can be classified as either:
\begin{itemize}
\item Geometric - based solely on the coordinate information of the vertices; 
\item Combinatorial - which attempt to group together highly connected vertices; 
\item Spectral - formulate the problem as the optimization of a discrete
quadratic function. The relaxed counterpart of the discrete problem
becomes a continuous one, which can be solved by computing the second
eigenvector of the discrete Laplacian of the graph;
\item Multilevel methods - a sequence of smaller graphs is constructed in
order to produce a similar coarser graph. The initial bisection is
performed on the smallest of these graphs. Finally, the graph is uncoarse
and partition refinement is performed on each of the coarse graph.
\end{itemize}
Each method has its advantages and disadvantages, and many of them
are described in \cite{Karypis1}, where we can find a detailed description
of several different methods in each of these classes. Combining those
methods is a common strategy to overcome the disadvantages. For instance,
spectral schemes can use eigenvectors to produce coordinate information
for vertices. Geometric methods can then use these coordinates to
partition the graph. Usually, for each application it is unclear which
method is better. There are many factors to be considered: degree
of parallelism, run time, quality of the cut produced. In \cite{karapis},
the author evaluate different aspects for many combinations of methods.
In general, it is agreed that spectral methods are good, specially
multilevel spectral bisection. 

In this paper, we aggregate more information present in the spectra
to improve the tradition spectral bisection algorithm (SB) and produce
a new graph bisection algorithm. While SB makes use of one eigenvector
only, the new algorithm uses two eigenvectors, which allows us to
returns a partition with cut size smaller or equal to the SB cut size.
Besides, the additional running time of computing an extra eigenvector
is rather small compared to the overall running time of SB. 

One one hand, we are specially concerned with the theoretical relations
of eigenvectors and cuts on graphs, and also show there still more
to be understood about these relations. Therefore, we do not intend
to make an extensive comparison between different classes of algorithms
and the new one, since the new algorithm is guaranteed to return a
cut that is not worse than the one of SB, at a cost of a rather small
running time. Nevertheless, we present some numerical results comparing
the quality of the cut between the new algorithm and SB. It is worth
it to mention that there is no restriction on using the new algorithm
in combination with other methods, and we expect that the new algorithm
improves the existing mixed methods that make use of the traditional
SB.

To reach our goal, we investigate properties of certain configurations
of a graph partition which are related to the minimality of a cut
and the structure of the graph, and we prove several results on that.
Such configurations, that we call \textit{organized partitions}, are
shown in this paper to be related to eigenvectors of the Laplacian
matrix. Turns out that organized partitions are relate to the maximum
cut problem as well, as we will show in section \ref{sec:Organized}.

Finally, we combine the organized partition, the third, and the second
eigenvector to construct an algorithm that approximates a minimum
graph bisection. For this algorithm, it is proven that the resulting
partition has number of edges smaller or equal than the classical
spectral bisection algorithm. Besides, we provide a second algorithm
that can produce a smaller cut, given a known cut, a procedure known
as refining a partition. There are several multilevel algorithms \cite{BJ2,Cheng,Karapis0,Hendrickson1}
that further refine the partition during the uncoarsening phase. The
second algorithm presented in this paper refines the partition by
making use the information about the organized partition present in
the third eigenvector.

The rest of the paper is organized as follow: properties of organized
partitions are investigated on section \ref{sec:Organized} and related
to minimum and maximum cuts on graphs. In section \ref{sec:IntegerProgramFormulation},
we connect organized partitions with spectral properties of graphs,
and we prove bounds on the minimum cut in terms of these properties.
In section \ref{sec:algorithm}, we derive both algorithms, the first
improving SB, and the second producing a smaller cut based on a given
one. In section \ref{sec:Experimentalresults}, we present some experimental
results comparing the quality of partitions returned by SB and the
new algorithm.

\section{\label{sec:Organized}Organized Partitions}

\begin{figure}[H]
\includegraphics{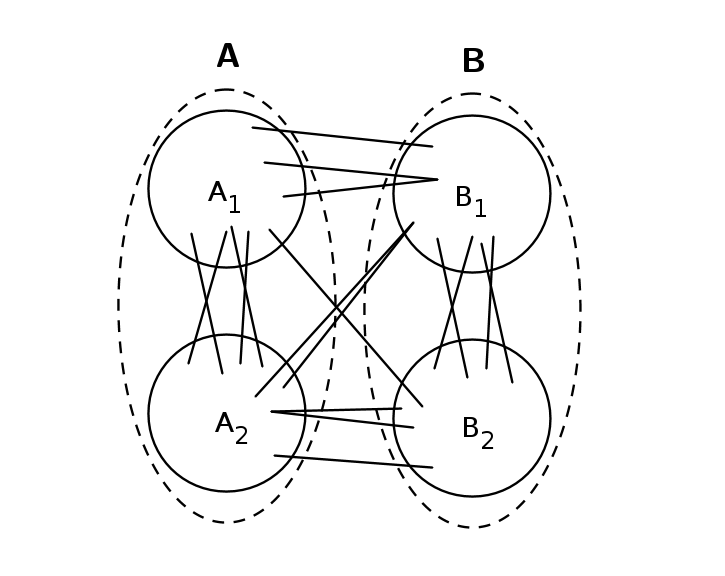}

\caption{Organized partition}
\label{fig:organized}
\end{figure}

Let $G=(V,E)$ be a connected graph with $4n$ vertices. Consider
a cut $\left\{ A,B\right\} $ of the vertex set $V$ such that $\left|A\right|=\left|B\right|$.
Such cut is also known as a balanced cut or a bisection. In this paper
we deal only with balanced cuts, thus from now on we will simply refer
to it simply as a cut. Let $A=A_{1}\cup A_{2}$ and $B=B_{1}\cup B_{2}$.
Now create a new partition of vertices $\mathcal{C}=\{A_{1},A_{2},B_{1},B_{2}\}$.
Here $E(X,Y)$ denotes the number of edges between the set of vertices
$A$ and $B$. We say that the partition $\mathcal{C}$ is organized
whenever 
\begin{equation}
E(A_{1},A_{2})+E(B_{1},B_{2})-E(A_{1},B_{1})-E(A_{2},B_{2})\label{eq:def1}
\end{equation}
is minimum among all subsets with $\left|A_{1}\right|=\left|A_{2}\right|=\left|B_{1}\right|=\left|B_{2}\right|$.
See Figure \ref{fig:organized}, which depicts the partition in question. 

It is worth mentioning that saying $\mathcal{C}$ is organized is
equivalent to say that 
\begin{equation}
E(A_{1},B_{2})+E(A_{2},B_{1})+E(A_{1},A_{2})+E(B_{1},B_{2})\label{eq:def2}
\end{equation}
is minimum among the prescribed sets. To see that, we notice that
\begin{eqnarray*}
 &  & E(A_{1},B_{2})+E(A_{2},B_{1})+E(A_{1},A_{2})+E(B_{1},B_{2})\\
 &  & =E(A,B)+E(A_{1},A_{2})+E(B_{1},B_{2})-E(A_{1},B_{1})-E(A_{2},B_{2}).
\end{eqnarray*}
Since $A$ and $B$ are fixed, we know that $E(A,B)$ is fixed too.
Thus the same subsets that minimize (\ref{eq:def1}) also minimize
(\ref{eq:def2}). We will show later how organized partitions relate
to minimum and maximum cuts of graphs.

In this paper, we tacitly assume that any partition $\mathcal{C}$
has $\left|A\right|=\left|B\right|$ and $\left|A_{1}\right|=\left|A_{2}\right|=\left|B_{1}\right|=\left|B_{2}\right|$.
Now, given a cut $\{A,B\}$ we can compute the quantity
\[
D_{\mathcal{C}}=\min_{\substack{A=\bar{A_{1}}\cup\bar{A_{2}}\\
B=\bar{B_{1}}\cup\bar{B_{2}}\\
\left|\bar{A_{1}}\right|=\left|\bar{A_{2}}\right|\\
\left|\bar{B_{1}}\right|=\left|\bar{B_{2}}\right|
}
}E(\bar{A_{1}},\bar{A_{2}})+E(\bar{B_{1}},\bar{B_{2}})-E(\bar{A_{1}},\bar{B_{1}})-E(\bar{A_{2}},\bar{B_{2}}).
\]
In this notation the solution of the optimization problem $\mathcal{C}=\left\{ A_{1},A_{2},B_{1},B_{2}\right\} $
is an organized partition for $\{A,B\}$.

We say that $\mathcal{C}$ is a minimum cut whenever $E(A,B)$ is
minimum among all choices of $A$ and $B$ with $\left|A\right|=\left|B\right|$.
In section \ref{sec:IntegerProgramFormulation}, we will see that
the quantity $D_{\mathcal{C}}$ relates with the eigenvalues of the
Laplacian matrix whenever $\mathcal{C}$ is a minimum cut.

The next Theorem provides a necessary condition for a cut to be minimum
or maximum from the perspective of its organized partition.
\begin{thm}
\label{thm:NotMin}Let $\left\{ A,B\right\} $ be any cut with organized
partition $\mathcal{C}$. If $D_{\mathcal{C}}<0$, then $\left\{ A,B\right\} $
is not a minimum cut. If $D_{\mathcal{C}}>0$, then $\left\{ A,B\right\} $
is not a maximum cut.
\end{thm}
\begin{proof}
Let $\mathcal{C}=\left\{ A_{1},A_{2},B_{1},B_{2}\right\} $. Notice
that 
\begin{eqnarray}
E(A,B) & = & E(A_{1},B_{1})+E(A_{1},B_{2})+E(A_{2},B_{1})+E(A_{2},B_{2}).\label{eq:ineq1}
\end{eqnarray}
Besides, 
\begin{equation}
E(A_{1}\cup B_{1},A_{2}\cup B_{2})=E(A_{1},A_{2})+E(B_{1},B_{2})+E(A_{2},B_{1})+E(A_{2},B_{2}).\label{eq:ineq2}
\end{equation}
Now, if $D_{\mathcal{C}}<0$ and from its the definition, we have
\[
E(A_{1},A_{2})+E(B_{1},B_{2})<E(A_{1},B_{1})+E(A_{2},B_{2}).
\]
This together with (\ref{eq:ineq1}) and (\ref{eq:ineq2}), gives
us 
\begin{eqnarray*}
E(A_{1}\cup B_{1},A_{2}\cup B_{2}) & < & E(A_{1},B_{1})+E(A_{2},B_{2})+E(A_{2},B_{1})+E(A_{2},B_{2})\\
 & = & E(A,B).
\end{eqnarray*}
Therefore, $\left\{ A,B\right\} $ is not a minimum cut. 

If $D_{\mathcal{C}}>0$, then 
\[
E(A_{1},A_{2})+E(B_{1},B_{2})>E(A_{1},B_{1})+E(A_{2},B_{2}).
\]
Similarly as before, that gives us
\[
E(A_{1}\cup B_{1},A_{2}\cup B_{2})>E(A,B).
\]
Thus, $\left\{ A,B\right\} $ is not a maximum cut. That finishes
the proof.
\end{proof}
In fact, the proof reveals a way to construct a better cut. That is
one of the fundamental ideas behind the algorithm we provide in section
\ref{sec:algorithm}. We explicit this construction in the form of
Corollary.
\begin{cor}
\label{cor:improveCut}If a cut $\left\{ A,B\right\} $ has $D_{\mathcal{C}}<0$,
then $E(A_{1}\cup B_{1},A_{2}\cup B_{2})<E(A,B)$. If $D_{\mathcal{C}}>0$,
then $E(A_{1}\cup B_{1},A_{2}\cup B_{2})>E(A,B)$.
\end{cor}
The next result gives some insights on how the organized partition
of a minimum/maximum cut looks like.
\begin{thm}
If $\left\{ A,B\right\} $ is a minimum cut, then its organized partition
satisfies $E(A_{1},A{}_{2})+E(B_{1},B_{2})\neq0$. If $\left\{ A,B\right\} $
is a maximum cut, then its organized partition satisfies $E(A_{1},B_{1})+E(A_{2},B_{2})\neq0$.
\end{thm}
\begin{proof}
Let $\left\{ A,B\right\} $ be a minimum cut and assume by contradiction
that $E(A_{1},A_{2})+E(B_{1},B_{2})=0$. We can assume that $E(A_{1},B_{1})+E(A_{2},B_{2})\neq0$,
otherwise the graph would be disconnected. Thus
\[
D_{\mathcal{C}}=E(A_{1},A_{2})+E(B_{1},B_{2})-E(A_{1},B_{1})-E(A_{2},B_{2})<0.
\]
Therefore, Theorem \ref{thm:NotMin} implies that $\left\{ A,B\right\} $
is not a minimum cut, which is a contradiction. 

If $\left\{ A,B\right\} $ is a maximum cut, assume by contradiction
that $E(A_{1},B_{1})+E(A_{2},B_{2})=0$. If $E(A_{1},A_{2})+E(B_{1},B_{2})=0$,
then the graph would be disconnected. Thus $E(A_{1},A_{2})+E(B_{1},B_{2})\neq0$,
and that gives us
\[
D_{\mathcal{C}}=E(A_{1},A_{2})+E(B_{1},B_{2})-E(A_{1},B_{1})-E(A_{2},B_{2})>0.
\]
Finally, Theorem \ref{thm:NotMin} implies that $\left\{ A,B\right\} $
is not a maximum cut, which is a contradiction. That finishes the
proof.
\end{proof}
Organized partitions also indicate conditions for which a graph has
more than one minimum or maximum cut and, if that is the case, how
to construct them.
\begin{thm}
Let $\left\{ A,B\right\} $ be any cut with organized partition $\mathcal{C=}\left\{ A_{1},A_{2},B_{1},B_{2}\right\} $.
If $D_{\mathcal{C}}=0$, then $E(A,B)=E(A_{1}\cup B_{1},A_{2}\cup B_{2})$. 
\end{thm}
\begin{proof}
From the definition of $D_{\mathcal{C}}$, we have 
\[
E(A_{1},A_{2})+E(B_{1},B_{2})=E(A_{1},B_{1})+E(A_{2},B_{2}).
\]
Thus, we can write
\begin{eqnarray*}
E(A_{1}\cup B_{1},A_{2}\cup B_{2}) & = & E(A_{1},B_{1})+E(A_{2},B_{2})+E(A_{2},B_{1})+E(A_{2},B_{2})\\
 & = & E(A,B).
\end{eqnarray*}
That finishes the proof.
\end{proof}
\begin{cor}
Let $\left\{ A,B\right\} $ be a minimum or a maximum cut. If $D_{\mathcal{C}}=0$,
then it is not unique.
\end{cor}
Thus, in some cases finding a organized partition can be useful to
construct a different minimum bisection whenever it is not unique.
On the other hand, for a graph with a unique minimum bisection, the
organized partition can be used to bound the size of the second minimum
bisection. As the next Theorem shows, a second minimum bisection is
not too far from the minimum whenever $D_{\mathcal{C}}$ is small.
\begin{thm}
Let $\left\{ A,B\right\} $ be a unique minimum cut and $\mathcal{C}=\{A_{1},A_{2},B_{1},B_{2}\}$
its organized partition. Let $\left\{ R,S\right\} $ be a second minimum
cut. Then 
\[
E(R,S)-E(A,B)\leq D_{\mathcal{C}}+1.
\]
\end{thm}
\begin{proof}
By Theorem \ref{thm:NotMin}, $D_{\mathcal{C}}\geq0$. Form a new
graph $G^{*}$ by adding $D_{\mathcal{C}}+1$ edges between $A_{1}$
and $B_{1}$. For this new graph it still holds that $\mathcal{C}=\{A_{1},A_{2},B_{1},B_{2}\}$
is an organized partition. Similarly, denoting by $D_{\mathcal{C}}^{*}$
and $E^{*}(A,B)$ the corresponding quantities in the graph $G^{*}$,
it holds that $D_{\mathcal{C}}^{*}=-1$ and $E^{*}(A,B)=E(A,B)+D_{\mathcal{C}}+1$. 

Now, assume by contradiction that $E(R,S)>D_{\mathcal{C}}+E(A,B)+1$.
Then we have $E^{*}(R,S)\geq E(R,S)>E^{*}(A,B).$ If we consider any
cut $\left\{ X,Y\right\} $ different than $\left\{ A,B\right\} $,
it holds that 
\[
E^{*}(X,Y)\geq E(X,Y)\geq E(R,S)>E^{*}(A,B),
\]
since $\left\{ R,S\right\} $ is a second minimum cut. This implies
that $\left\{ A,B\right\} $ is a minimum cut for $G^{*}$ as well.
By Theorem \ref{thm:NotMin}, this minimum cut satisfy $D_{\mathcal{C}}^{*}\geq0$.
That is a contradiction with $D_{\mathcal{C}}^{*}=-1$. Therefore,
$E(R,S)\leq D_{\mathcal{C}}+E(A,B)+1$, which finishes the proof.
\end{proof}

\section{\label{sec:IntegerProgramFormulation}Integer Program Formulation}

This section is dedicated to relate organized partitions with spectral
properties of the graph. We prove bounds on the minimum cut in terms
of these properties. In the next Theorem, we show how to construct
the organized partition of given cut. Turns out it suffices to solve
an integer program in terms of the Laplacian matrix of the graph.
\begin{thm}
\label{thm:ProgOrganized}Let $\left\{ A,B\right\} $ be any bisection
of a graph $G$ and denote by $y$ be the vector with entries 
\[
y_{i}=\begin{cases}
1/\sqrt{n} & \text{ if }i\in A\\
-1/\sqrt{n} & \text{ if }i\in B.
\end{cases}
\]
Let $L$ be the Laplacian matrix of the $G$. Then 
\begin{equation}
\frac{4}{n}\left(E(A,B)+D_{\mathcal{C}}\right)=\min_{\substack{x^{T}\mathbf{1}=0\\
\left\Vert x\right\Vert =1\\
y^{T}x=0\\
x_{i}\in\left\{ 1/\sqrt{n},-1/\sqrt{n}\right\} 
}
}x^{T}Lx.\label{eq:IntegerProgl3}
\end{equation}
Furthermore, each solution $\bar{x}$ of (\ref{eq:IntegerProgl3})
prescribes an organized partition for $\left\{ A,B\right\} $ as follow

\[
\bar{x_{i}}=\begin{cases}
1/\sqrt{n} & \text{ }i\in A_{1}\cup B_{1}\\
-1/\sqrt{n} & \text{ }i\in A_{2}\cup B_{2}.
\end{cases}
\]
\end{thm}
\begin{proof}
Let $A_{1}$,$A_{2}$, $B_{1}$and $B_{2}$ be disjoint sets such
that $A_{1}\cup A_{2}=A$ and $B_{1}\cup B_{2}=B$, with $\left|A_{1}\right|=\left|A_{2}\right|$
and $\left|B_{1}\right|=\left|B_{2}\right|$. Define the vector $x$
with entries 
\[
x_{i}=\begin{cases}
1/\sqrt{n} & i\in A_{1}\cup B_{1}\\
-1/\sqrt{n} & i\in A_{2}\cup B_{2}
\end{cases}.
\]
Clearly $x^{T}\mathbf{1}=0$, $\left\Vert x\right\Vert =1$, and $y^{T}x=0.$ 

Now, we can write $x^{T}Lx$ in terms of the partition $\left\{ A_{1},A_{2},B_{1},B_{2}\right\} $
as 
\begin{eqnarray*}
x^{T}Lx & = & \sum_{(i,j)\in E}\left(x_{i}-x_{j}\right)^{2}\\
 & = & \sum_{\substack{(i,j)\in E\\
i\in A_{1},j\in B_{1}
}
}\left(x_{i}-x_{j}\right)^{2}+\sum_{\substack{(i,j)\in E\\
i\in A_{2},j\in B_{2}
}
}\left(x_{i}-x_{j}\right)^{2}+\sum_{\substack{(i,j)\in E\\
i\in A_{1},j\in B_{2}
}
}\left(x_{i}-x_{j}\right)^{2}\\
 &  & +\sum_{\substack{(i,j)\in E\\
i\in A_{2},j\in B_{2}
}
}\left(x_{i}-x_{j}\right)^{2}+\sum_{\substack{(i,j)\in E\\
i\in A_{1},j\in A_{2}
}
}\left(x_{i}-x_{j}\right)^{2}+\sum_{\substack{(i,j)\in E\\
i\in B_{1},j\in B_{2}
}
}\left(x_{i}-x_{j}\right)^{2}\\
 & = & \frac{4}{n}\left(E(A_{1},B_{2})+E(A_{2},B_{1})+E(A_{1},A_{2})+E(B_{1},B_{2})\right),
\end{eqnarray*}
since the first two sums are zero. That gives us

\[
\frac{n}{4}x^{T}Lx=E(A,B)+E(A_{1},A_{2})+E(B_{1},B_{2})-E(A_{1},B_{1})-E(A_{2},B_{2}),
\]
for each choice of partition $\left\{ A_{1},A_{2},B_{1},B_{2}\right\} $.
Therefore, in view of the definition of $D_{\mathcal{C}}$, we have

\[
\min_{\substack{x^{T}\mathbf{1}=0\\
\left\Vert x\right\Vert =1\\
y^{T}x=0\\
x_{i}\in\left\{ 1/\sqrt{n},-1/\sqrt{n}\right\} 
}
}x^{T}Lx=\frac{4}{n}\left(E(A,B)+D_{\mathcal{C}}\right).
\]

Besides, by the construction of the feasible set of solutions, $\bar{x}$
indicates the organized partition of $\left\{ A,B\right\} $. That
finishes the proof.
\end{proof}
Thus, whenever $\left\{ A,B\right\} $ is a minimum cut, the minimum
of (\ref{eq:IntegerProgl3}) reduces to $\frac{4}{n}\left(MinCut(G)+D_{\mathcal{C}}\right)$.

In the work of \cite{DonathHoffman2} the authors proved the inequality
\begin{equation}
MinCut(G)\geq\frac{n}{4}\lambda_{2}.\label{eq:boundAlg}
\end{equation}
In light of the concept of organized partitions we can go further
on the relation between minimum cuts and eigenvalues of the Laplacian
matrix and prove the next result.
\begin{thm}
\label{thm:minL2L3}Let $\mathcal{C}$ be a organized partition of
a minimum cut. Then 

\[
MinCut(G)\geq\frac{n}{8}\left(\lambda_{2}+\lambda_{3}\right)-\frac{D_{\mathcal{C}}}{2}.
\]
\end{thm}
\begin{proof}
Define the vector $y$ with entries 

\begin{equation}
y_{i}=\begin{cases}
1/\sqrt{n} & i\in A\\
-1/\sqrt{n} & i\in B
\end{cases}\label{eq:vecy}
\end{equation}
Clearly $y^{T}\mathbf{1}=0$ and $\left\Vert y\right\Vert =1$. Thus,
we can write
\[
y^{T}Ly=\sum_{(i,j)\in E}\left(y_{i}-y_{j}\right)^{2}=\sum_{\substack{(i,j)\in E\\
i\in A,j\in B
}
}\left(y_{i}-y_{j}\right)^{2}+\sum_{\substack{(i,j)\in E\\
i\in A,j\in A
}
}\left(y_{i}-y_{j}\right)^{2}+\sum_{\substack{(i,j)\in E\\
i\in B,j\in B
}
}\left(y_{i}-y_{j}\right)^{2}.
\]
Notice the sum over the edges with both endpoints in the same set
vanishes. Thus, we have
\[
y^{T}Ly=\sum_{\substack{(i,j)\in E\\
i\in A,j\in B
}
}\left(1/\sqrt{n}-(-1/\sqrt{n})\right)^{2}=\frac{4}{n}E(A,B).
\]
An important idea here is that a minimum cut is achieved if we take
the minimum over all prescribed vectors, i.e., 
\begin{equation}
MinCut(G)=\frac{n}{4}\min_{\substack{y^{T}\mathbf{1}=0,\left\Vert y\right\Vert =1\\
y_{i}\in\left\{ 1/\sqrt{n},-1/\sqrt{n}\right\} 
}
}y^{T}Ly.\label{eq:prog1}
\end{equation}

Now, we apply Theorem \ref{thm:ProgOrganized} for the vector $\bar{y}$
that solves the minimization problem (\ref{eq:prog1}). Thus, we can
solve the sum of minimization problems as

\[
\min_{\substack{y^{T}\mathbf{1}=0,\left\Vert y\right\Vert =1\\
y_{i}\in\left\{ 1/\sqrt{n},-1/\sqrt{n}\right\} 
}
}y^{T}Ly+\min_{\substack{x^{T}\mathbf{1}=0\\
\left\Vert x\right\Vert =1\\
\bar{y}^{T}x=0\\
x_{i}\in\left\{ 1/\sqrt{n},-1/\sqrt{n}\right\} 
}
}x^{T}Lx=\frac{4}{n}\left(2MinCut(G)+D_{\mathcal{C}}\right).
\]
Equivalently, we can write 
\[
MinCut(G)=\frac{n}{8}\min_{\substack{y^{T}\mathbf{1}=0,\left\Vert y\right\Vert =1\\
y_{i}\in\left\{ 1/\sqrt{n},-1/\sqrt{n}\right\} 
}
}y^{T}Ly+\min_{\substack{x^{T}\mathbf{1}=0\\
\left\Vert x\right\Vert =1\\
\bar{y}^{T}x=0\\
x_{i}\in\left\{ 1/\sqrt{n},-1/\sqrt{n}\right\} 
}
}x^{T}Lx-\frac{D_{\mathcal{C}}}{2}.
\]
Thus, if we drop the constraint $y_{i},x_{i}\in\left\{ 1/\sqrt{n},-1/\sqrt{n}\right\} $
and consider all $x,y\in\mathbb{R}^{n}$, we find the inequality
\begin{eqnarray*}
MinCut(G) & \geq & \frac{n}{8}\min_{\substack{y^{T}\mathbf{1}=0\\
\left\Vert y\right\Vert =1
}
}y^{T}Ly+\min_{\substack{x^{T}\mathbf{1}=0\\
\left\Vert x\right\Vert =1\\
\bar{y}^{T}x=0
}
}x^{T}Lx-\frac{D_{\mathcal{C}}}{2}\\
 & = & \frac{n}{8}\min_{\substack{y^{T}\mathbf{1}=x^{T}\mathbf{1}=0\\
\left\Vert y\right\Vert =\left\Vert x\right\Vert =1\\
y^{T}x=0
}
}y^{T}Ly+x^{T}Lx-\frac{D_{\mathcal{C}}}{2}\\
 & = & \frac{n}{8}\left(\lambda_{2}+\lambda_{3}\right)-\frac{D_{\mathcal{C}}}{2}.
\end{eqnarray*}
That finishes the proof.
\end{proof}
It is worth mentioning that if we drop the constraint $y_{i}\in\left\{ 1/\sqrt{n},-1/\sqrt{n}\right\} $
in the minimization problem (\ref{eq:prog1}), then we precisely obtain
the lower bound (\ref{eq:boundAlg}) as the authors in \cite{DonathHoffman2}.
Besides, whenever 
\[
D_{\mathcal{C}}<\frac{n}{4}\left(\lambda_{3}-\lambda_{2}\right),
\]
Theorem \ref{thm:minL2L3} provides a tighter lower bound on the minimum
cut of a graph. Intuitively, it means that an optimization problem
that considers both $\lambda_{2}$ and $\lambda_{3}$ is more likely
to reveal a minimum cut than a problem that considers only $\lambda_{2}$.

We finish this section with a result that summarizes all its underlying
ideas.
\begin{thm}
\label{thm:integerFULL}For a graph with Laplacian matrix $L$, the
solution $\left(\bar{x},\bar{y}\right)$ of the problem 
\[
\min_{\substack{y^{T}\mathbf{1}=0,\left\Vert y\right\Vert =1\\
y_{i}\in\left\{ 1/\sqrt{n},-1/\sqrt{n}\right\} 
}
}y^{T}Ly+\min_{\substack{x^{T}\mathbf{1}=0\\
\left\Vert x\right\Vert =1\\
\bar{y}^{T}x=0\\
x_{i}\in\left\{ 1/\sqrt{n},-1/\sqrt{n}\right\} 
}
}x^{T}Lx
\]
constructs a minimum cut $\left\{ A,B\right\} $ together with its
organized partition\\
 $\mathcal{C}=\{A_{1},A_{2},B_{1},B_{2}\}$, as follow: 
\[
\bar{y_{i}}=\begin{cases}
1/\sqrt{n} & i\in A\\
-1/\sqrt{n} & i\in B
\end{cases}\text{ and }\bar{x_{i}}=\begin{cases}
1/\sqrt{n} & \text{ if }i\in A_{1}\cup B_{1}\\
-1/\sqrt{n} & \text{ if }i\in A_{2}\cup B_{2}.
\end{cases}
\]
\end{thm}
\begin{proof}
Follows from equation (\ref{eq:prog1}) and Theorem \ref{thm:ProgOrganized}.
\end{proof}

\section{\label{sec:algorithm}Derivation of the Algorithms}

In this section we provide an intuitive description of the main ideas
behind our new algorithm, which turns out to arise from the theoretical
background developed in the previous sections. We do that by showing
how to improve the bisection provided by the traditional SB algorithm
by means of properties of organized partitions. We will prove that
there are infinite many solutions for the minimization problem that
finds the organized partition of a cut, if we apply relaxation. Thus,
these solutions constructs better candidates for a minimum cut. First,
we consider some examples where SB fails to approximate a good bisection.

\begin{figure}[H]
\includegraphics[scale=0.4]{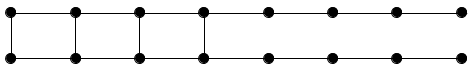}

\caption{Roach graph on 16 vertices}
\label{fig:Roach}
\end{figure}

As an approximation algorithm, SB sometimes provides a cut that is
too far from optimal. There are investigations about this phenomenon,
and the best known example where SB fails is given by the roach graph,
due to Guattery and Miller \cite{Guattery}. The roach graph consists
of two path graphs with the same even size connected by a few edges,
as illustrated in Figure \ref{fig:Roach}. 

This is a very good example which seems to be taylor made to defeat
SB. The roach graph is an important example not only because SB provides
a cut that is far from optimal, in fact it is the prototype of many
cases where this algorithm gives a very bad result. Let us look closer
to what is happening with the algorithm on this kind of graph.

\begin{figure}[H]
\includegraphics[scale=0.4]{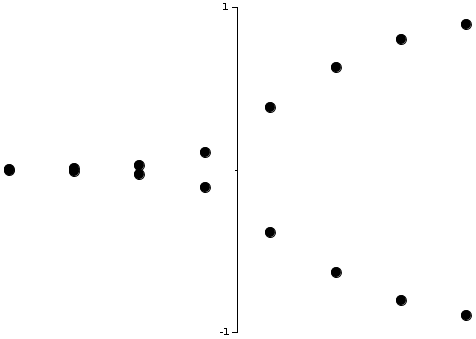}

\caption{The $y$ eigenvector for the roach graph}
\label{fig:roachY}
\end{figure}

For a roach graph the minimum bisection consists of two edges separating
the antennae - the pending paths on the right side of Figure \ref{fig:Roach}.
But that is not what SB returns. Taking a roach graph on 16 vertices,
we label the upper and lower path from 1 to 8 and 9 to 16, respectively.
For this ordering, its eigenvector associated with $\lambda_{2}$
is approximately given by{\small{}
\begin{eqnarray*}
y & = & [-0.0028-0.0083-0.0295-0.1068-0.3869796-0.6270-0.8024-0.8948\\
 &  & 0.0028\text{ }0.0083\text{ }0.0295\text{ }0.1068\text{ }0.3869\text{ }0.6270\text{ }0.8024\text{ }0.8948]^{T}.
\end{eqnarray*}
}Now, we can plot the entries of $y$ displayed in Figure \ref{fig:roachY}.
The upper path corresponds to the points above the origin and the
lower path bellow it. SB will split the graph in two paths, which
provides a cut with 4 edges, which is not a maximum cut. In \cite{Guattery}
the authors showed this is true for the whole class of roach graphs,
therefore showing a class of graphs where the resulting bisection
from SB is far from optimal, i.e., with a bisection of order $\mathcal{O}(n)$.

This is the prototype of what happens with SB when it returns a wrong
bisection. In view of this problem, it is natural to ask how to overcome
this pathology on the SB algorithm. Here we show how that can be done
using the concept of organized partitions.

In light of Corollary \ref{cor:improveCut} if a cut has $D_{\mathcal{C}}<0$,
then its organized partition can be used to construct a smaller cut.
Thus, it would be useful to have an algorithm that approximates a
minimum cut and which computes its organized partition as well. That
means if we could solve both problems simultaneously, then we can
obtain a better cut than the original algorithm, whenever this cut
has $D_{\mathcal{C}}<0$. 

That is the case of the roach graph and many examples of this nature.
Notice that the cut provided by SB for the roach graph on 16 vertices
is $\mathcal{C}=\left\{ A_{1},A_{2},B_{1},B_{2}\right\} $, where\\
 $A_{1}=\left\{ v_{1},\ldots,v_{4}\right\} $, $A_{2}=\left\{ v_{5},\ldots,v_{8}\right\} $,
$B_{1}=\left\{ v_{9},\ldots,v_{12}\right\} $, and $B_{2}=\left\{ v_{13},\ldots,v_{16}\right\} $.
Therefore, we have
\[
D_{\mathcal{C}}=E(A_{1},A_{2})+E(B_{1},B_{2})-E(A_{1},B_{1})-E(A_{2},B_{2})=1+1-4-0,
\]
which gives us the desired property $D_{\mathcal{C}}<0$. For this
reason, the organized partition of this cut will provide a smaller
bisection.

Now, let us see what the eigenvector of $\lambda_{3}$ tells about
the organized partition. Theorem \ref{thm:integerFULL} constructs
the organized partition based on the solution of an integer program.
Theorem \ref{thm:minL2L3} and its proof indicate that the eigenvectors
of $\lambda_{2}$ and $\lambda_{3}$ can be used to approximate the
solution. Thus, if we drop the constraints on $x$ and $y$ putting
$x,y\in\mathbb{R}^{n}$, it is expected that the solution of the new
program 

\begin{equation}
\min_{\substack{y^{T}\mathbf{1}=x^{T}\mathbf{1}=0\\
\left\Vert y\right\Vert =\left\Vert x\right\Vert =1\\
y^{T}x=0
}
}y^{T}Ly+x^{T}Lx\label{eq:Prog-eigenvectors}
\end{equation}
approximates the minimum cut and its organized partition by the eigenvector
$x$ associated with $\lambda_{3}$.

For the same roach graph, that eigenvector is approximately 
\begin{eqnarray*}
x & = & [-0.6935\text{ }-0.5879\text{ }-0.3928\text{ }-0.1379\text{ }0.1379\text{ }0.3928\text{ }0.5879\text{ }0.6935\text{ }\\
 &  & -0.6935\text{ }-0.5879-0.3928\text{ }-0.1379\text{ }\text{ }0.1379497\text{ }\text{ }0.3928475\text{ }\text{ }0.5879\text{ }\text{ }0.6935]^{T}.
\end{eqnarray*}

Notice that if we use $x$ as an approximation for the integer solution
of the program in Theorem \ref{thm:integerFULL}, then $x$ induces
the correct organized partition\\
 $\mathcal{C}=\left\{ A_{1},A_{2},B_{1},B_{2}\right\} $ as described
above. Here we simply used the entries of $x$ as an approximation
for the integer solution 
\[
\bar{x_{i}}=\begin{cases}
1/\sqrt{n} & \text{ if }i\in A_{1}\cup B_{1}\\
-1/\sqrt{n} & \text{ if }i\in A_{2}\cup B_{2}.
\end{cases}
\]

Since $D_{\mathcal{C}}<0$, this implies that we can construct a smaller
bisection than the one provided by SB by using the eigenvector $x$.
More precise, by Corollary \ref{cor:improveCut} the partition $\left\{ A_{1}\cup B_{1},A_{2}\cup B_{2}\right\} $
gives a smaller bisection. In fact, this is the minimum bisection
for the roach graph.

\begin{figure}[H]
\includegraphics[scale=0.4]{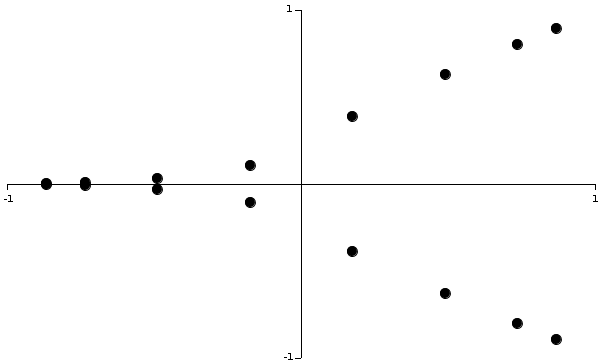}

\caption{Roach graph with $x$ and $y$ as embedding}
\label{fig:roachXandY}
\end{figure}

Figure \ref{fig:roachXandY} depicts the underlying idea behind the
proof of Theorem \ref{thm:minL2L3}. We plotted points using the entries
of both eigenvectors of the roach graph $x$ and $y$ as coordinates.
There, each point corresponds to a vertex. It is clear to see that
if we separate the vertices by the signs of the coordinates in $x$,
then we would get the minimum cut. 

The previous discussion suggests to consider both eigenvectors in
a new algorithm, in the sense either $x$ or $y$ will approximate
a minimum cut. Essentially, when $y$ gives a cut with $D_{\mathcal{C}}<0$,
we can appeal to the cut provided by $x$. Thus it would suffices
to check which one gives a better cut. Actually, this neat idea can
be taken further when we look from the perspective of integer programming. 

As we will see in the next Theorem, certain specific linear combinations
of $x$ and $y$ are solutions for (\ref{eq:Prog-eigenvectors}) as
well. Thus, those new solutions can be used to approximate a minimum
bisection.
\begin{thm}
\label{thm:rotationOfSolutions}Let $x$ and $y$ be a solution of
(\ref{eq:Prog-eigenvectors}). Let $\theta\in[0,2\pi)$ and let $u=cos\theta x+sin\theta y$
and $v=sin\theta x-cos\theta y$. Then $u$ and $v$ is a solution
of (\ref{eq:Prog-eigenvectors}). 
\end{thm}
\begin{proof}
We proceed by showing that $x^{T}Lx+y^{T}Ly=u^{T}Lu+v^{T}Lv$. Hence,
we write 
\begin{eqnarray*}
u^{T}Lu & = & (\mbox{cos}\theta x+\mbox{sin}\theta y)^{T}L(\mbox{cos}\theta x+\mbox{sin}\theta y)\\
 & = & \mbox{cos}\theta x^{T}L\mbox{cos}\theta x+\mbox{sin}\theta y^{T}L\mbox{sin}\theta y+2\mbox{sin}\theta y^{T}L\mbox{cos}\theta x.
\end{eqnarray*}
Also, we can write

\begin{eqnarray*}
v^{T}Lv & = & (\mbox{sin}\theta x-\mbox{cos}\theta y)^{T}L(\mbox{sin}\theta x-\mbox{cos}\theta y)\\
 & = & \mbox{sin}\theta x^{T}L\mbox{sin}\theta x+\mbox{cos}\theta y^{T}L\mbox{cos}\theta y-2\mbox{sin}\theta y^{T}L\mbox{cos}\theta x.
\end{eqnarray*}
Therefore, we obtain 
\begin{eqnarray*}
u^{T}Lu+v^{T}Lv & = & \mbox{cos}\theta x^{T}L\mbox{cos}\theta x+\mbox{sin}\theta y^{T}L\mbox{sin}\theta y+2\mbox{sin}\theta y^{T}L\mbox{cos}\theta x\\
 &  & +\mbox{sin}\theta x^{T}L\mbox{sin}\theta x+\mbox{cos}\theta y^{T}L\mbox{cos}\theta y-2\mbox{sin}\theta y^{T}L\mbox{cos}\theta x\\
 & = & \mbox{co\ensuremath{s^{2}}}\theta x^{T}Lx+\mbox{si\ensuremath{n^{2}}}\theta y^{T}Ly+\mbox{si\ensuremath{n^{2}}}\theta x^{T}Lx+\mbox{co\ensuremath{s^{2}}}\theta y^{T}Ly\\
 & = & (\mbox{co\ensuremath{s^{2}}}\theta+\mbox{si\ensuremath{n^{2}}}\theta)(x^{T}Lx+y^{T}Ly)\\
 & = & x^{T}Lx+y^{T}Ly.
\end{eqnarray*}
It follows that $u$ and $v$ is also a minimizer of (\ref{eq:Prog-eigenvectors}). 

It remains to verify that $u$ and $v$ satisfy the constraints $u^{T}\mathbf{1}=0$,
$v^{T}\mathbf{1}=0$, $\left\Vert u\right\Vert =\left\Vert v\right\Vert =1$
and $u^{T}v=0$. 

To see that $u^{T}\mathbf{1}=0$, we notice that $u^{T}\mathbf{1}=\mbox{cos}\theta x^{T}\mathbf{1}+\mbox{sin}\theta y^{T}\mathbf{1}=0$.
Now, using the fact that $x^{T}x=y^{T}y=1$ and $y^{T}x=0$, we can
write
\begin{eqnarray*}
u^{T}u & = & (\mbox{cos}\theta x+\mbox{sin}\theta y)^{T}(\mbox{cos}\theta x+\mbox{sin}\theta y)\\
 & = & \mbox{co\ensuremath{s^{2}}}\theta x^{T}x+\mbox{si\ensuremath{n^{2}}}\theta y^{T}y+2\mbox{sin}\theta\mbox{cos}\theta y^{T}x\\
 & = & \mbox{co\ensuremath{s^{2}}}\theta+\mbox{si\ensuremath{n^{2}}}\theta=1,
\end{eqnarray*}
which implies $\left\Vert u\right\Vert =1$. Similarly, we obtain
$v^{T}\mathbf{1}=0$ and $\left\Vert v\right\Vert =1$.

To show $u^{T}v=0$, again we use the fact that $x^{T}x=y^{T}y=1$
and $y^{T}x=0$ 
\begin{eqnarray*}
u^{T}v & = & (\mbox{cos}\theta x+\mbox{sin}\theta y)^{T}(\mbox{sin}\theta x-\mbox{cos}\theta y)\\
 & = & \mbox{cos}\theta\mbox{sin}\theta x^{T}x-\mbox{co\ensuremath{s^{2}}}\theta x^{T}y+\mbox{si\ensuremath{n^{2}}}\theta x^{T}y-\mbox{sin}\theta\mbox{cos}\theta y^{T}y\\
 & = & \mbox{cos}\theta\mbox{sin}\theta x^{T}x-\mbox{sin}\theta\mbox{cos}\theta y^{T}y=0.
\end{eqnarray*}
This concludes the proof.
\end{proof}
By constructing a infinite set of solutions for the problem (\ref{eq:Prog-eigenvectors}),
the last Theorem introduces a degree of freedom in the solutions of
(\ref{eq:Prog-eigenvectors}). We can explore this degree of freedom
in order to create different bisections. As discussed before, solutions
of (\ref{eq:Prog-eigenvectors}) can be used to approximate a minimum
bisection and its organized partition. However, there are infinite
$u$ and $v$ described in the last Theorem. Naturally, all of them
can be used to approximate a minimum bisection. That is a key idea
in the algorithm presented next. The next Theorem show how to construct
different $n$ different bisections based on the solutions of solutions
of (\ref{eq:Prog-eigenvectors}).
\begin{thm}
Let $x$ and $y$ be solutions of (\ref{eq:Prog-eigenvectors}). For
each pair $x_{i}$ and $y_{i}$, $i=1,\ldots,n$, define the vector
$u=\frac{x_{i}}{\sqrt{x_{i}^{2}+y_{i}^{2}}}x+\frac{y_{i}}{\sqrt{x_{i}^{2}+y_{i}^{2}}}y$.
Then $u$ induces a bisection that approximates the vector $\bar{u_{i}}$
with entries 
\[
\bar{u_{i}}=\begin{cases}
1/\sqrt{n} & i\in A\\
-1/\sqrt{n} & i\in B
\end{cases}.
\]
\end{thm}
\begin{proof}
In order to construct different bisections using Theorem \ref{thm:rotationOfSolutions}
we need to choose $\theta\in[0,2\pi)$, then define $u$ and $v$,
and finally define a new partition $\{A,B\}$ based on $u$ and $v$.
To this end, consider the set of euclidean points $(x_{i},y_{i})$
given by the corresponding entries of the eigenvectors $x$ and $y$.
Choose a point $(x_{i},y_{i})$, and let $\theta_{i}$ be the angle
between the point $(x_{i},y_{i})$ and the abscissa. Now define $u$
and $v$ as in Theorem \ref{thm:rotationOfSolutions}, and let $(u_{i},v_{i})$
be points defined by the corresponding entries of $u$ and $v$. The
point $(u_{i},v_{i})$ is simply a rotation of angle $\theta_{i}$
for the point $(x_{i},y_{i})$. 

Now using the solution of (\ref{eq:Prog-eigenvectors}), we can approximate
the solution of the integer program in Theorem (\ref{thm:integerFULL}).
By Theorem (\ref{thm:integerFULL}), its solution defines a minimum
cut, and we can define the cut $\{A,B\}$ using the entries of $u$
as an approximation for 
\[
\bar{u_{i}}=\begin{cases}
1/\sqrt{n} & i\in A\\
-1/\sqrt{n} & i\in B
\end{cases}.
\]

Finally, to simplify the computation of $u$ we can calculate $\mbox{cos}\theta$
and $\mbox{sin}\theta$ instead of $\theta$. That follows straightforward
from 
\[
\mbox{cos}\theta=\frac{x_{i}}{\sqrt{x_{i}^{2}+y_{i}^{2}}}\text{ and }\mbox{sin}\theta=\frac{y_{i}}{\sqrt{x_{i}^{2}+y_{i}^{2}}}.
\]
That finishes the proof.
\end{proof}
Now we are ready to give the complete algorithm that approximates
a minimum bisection of a graph.

\begin{algorithm}[H]
\begin{algorithmic} 
\Require{G=(V,E)} 
\State Compute $y$ and $x$, the second and third smallest eigenvector of $L$.
\State Set $A$  with the $n/2$ vertices with largest $y_i$ and $B$ with the remaining vertices.
\For{$i=1,\ldots,n$}
\State $u=\frac{x_{i}}{\sqrt{x_{i}^{2}+y_{i}^{2}}}x+\frac{y_{i}}{\sqrt{x_{i}^{2}+y_{i}^{2}}}y$ 
\State Set $R$ with the $n/2$ vertices with largest $u_j$ and $S$ with the remaining vertices. 
\If{$E(R,S)<E(A,B)$} 
\State $A=R$  
\State $B=S$  
\EndIf 
\EndFor \\
\Return $\{A,B\}$ \\
\end{algorithmic}

\caption{Graph Bisection.}
\label{algo:algorithm1}
\end{algorithm}

\begin{figure}[H]
\subfloat[SB chooses vertices from a Fiedler vector]{\includegraphics[scale=0.45]{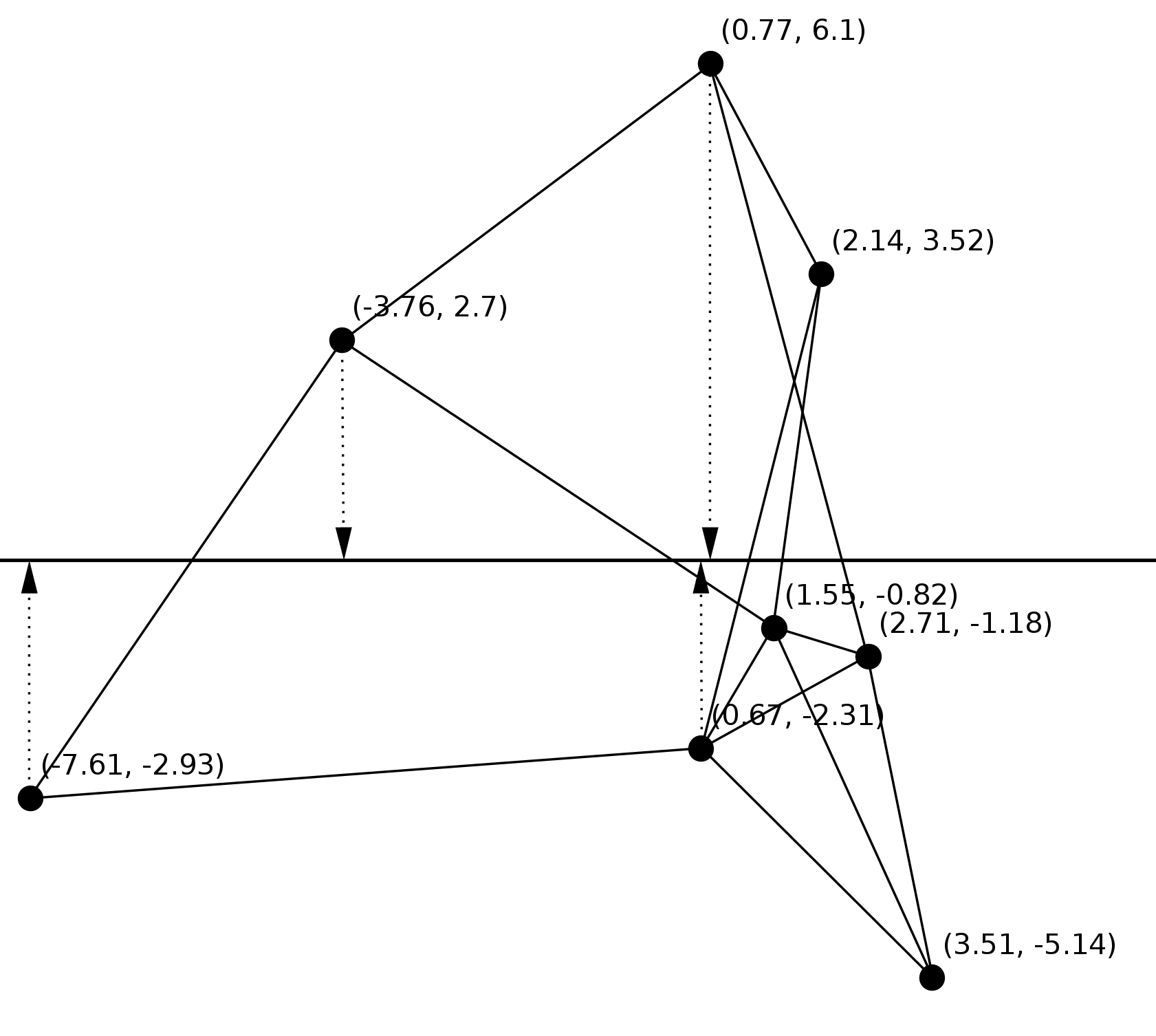}

\label{fig:exampleCutA}}\subfloat[Algorithm \ref{algo:algorithm1} considers both, the Fiedler vector
and the third eigenvector to choose vertices]{\includegraphics[scale=0.45]{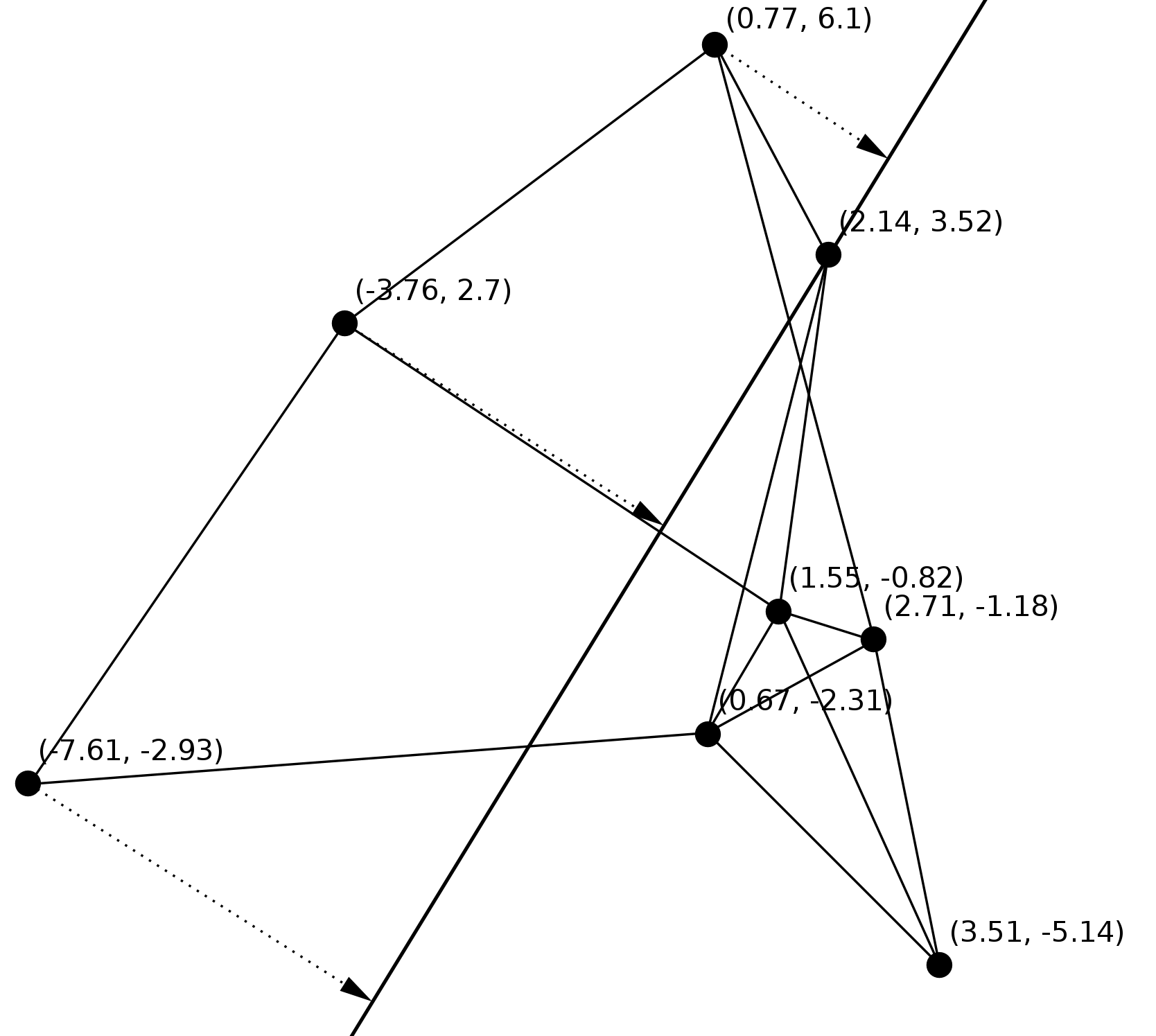}

\label{fig:exampleCutB}}

\caption{Different lines induce different bisections}
\end{figure}

As an illustration of Algorithm \ref{algo:algorithm1}, Figures \ref{fig:exampleCutA}
and \ref{fig:exampleCutB} show the same graph embedded on the coordinates
given by the second and the third eigenvalue. Figure \ref{fig:exampleCutA}
depicts the SB algorithm choosing a set of vertices based on a Fiedler
vector only. The straight line has the same direction of the Fiedler
vector. Since SB sorts the vertices based on the this vector and chooses
the top largest to construct the bisection, it is clear that it is
simply a projection of points along the straight line. As more linear
combinations of the Fiedler vector and the third eigenvector are considered,
different cuts are created. Figure \ref{fig:exampleCutB} depicts
the optimal choice of vertices induced by one of those linear combinations. 

Notice that the cut induced by $x$, the standard spectral bisection
solution, is among the possible cuts $\{R,S\}$ constructed by Algorithm
\ref{algo:algorithm1}. Therefore, the number of edges in the partition
provided by Algorithm \ref{algo:algorithm1} is not larger than the
one in the partition returned by SB, which leads us to the next Theorem.
\begin{thm}
The cut returned by Algorithm \ref{algo:algorithm1} has number of
edges smaller or equal than the number of edges in the SB partition.

$\square$
\end{thm}
For any roach graph its eigenvectors have the same shape of the previous
example with 16 vertices. That leads us to the next Theorem.
\begin{thm}
For any roach graph, Algorithm \ref{algo:algorithm1} returns a minimum
cut.
\end{thm}
\begin{proof}
By Lemma 5.1 of \cite{Guattery}, the third eigenvector of a roach
graph induces a cut separating the pending paths of the roach graph,
which is a minimum cut. This cut is among the possible cuts constructed
by Algorithm \ref{algo:algorithm1}. That finishes the proof.
\end{proof}
Now we will turn our attention to the derivation of an algorithm that
refines a given bisection. Since an organized partition can be used
to construct a better bisection, the next algorithm constructs an
approximation for an organized partition of a given bisection. In
the same fashion as in Algorithm \ref{algo:algorithm1}, these approximations
are candidates for a smaller cut.

Theorem \ref{thm:ProgOrganized}, provide us with a way to construct
the organized partition of a given cut. If $\left\{ A,B\right\} $
is the cut in question, we can denote by $y$ be the vector with entries
\[
y_{i}=\begin{cases}
1/\sqrt{n} & \text{ if }i\in A\\
-1/\sqrt{n} & \text{ if }i\in B.
\end{cases}
\]
Now, if we use relaxation on the set of solutions of the integer program
(\ref{eq:IntegerProgl3}) and drop the constraint $x_{i}\in\left\{ 1/\sqrt{n},-1/\sqrt{n}\right\} $,
we obtain the following program

\begin{equation}
\min_{\substack{x^{T}\mathbf{1}=0\\
\left\Vert x\right\Vert =1\\
y^{T}x=0
}
}x^{T}Lx.\label{eq:progImprove}
\end{equation}

The minimization problem (\ref{eq:progImprove}) is not an eigenvalue
problem anymore, because the vector $y$ is not necessarily an eigenvector
of the matrix $L$. However, it is easy to transform problem (\ref{eq:progImprove})
into a standard eigenvalue problem, as shown in \cite{Golub} by Gene
and Golub. Therefore, the solution of program (\ref{eq:progImprove})
can be used as an approximation for the organized partition: the half
largest entries of $x$ indicate the vertices in the set $A_{1}\cup B_{1}$
of the organized partition, and the other half indicates the remaining
vertices in the organized partition. Again, we can use linear combinations
of $x$ and $y$ to construct different approximations for the organized
partition. The algorithm can be described as follows.

\begin{algorithm}[H]
\begin{algorithmic} 
\Require{$G=(V,E),y$} 
\State Set $A$  with the $n/2$ vertices with largest $y_i$ and $B$ with the remaining vertices.
\State Compute $x$, the solution of 
$\min_{\substack{x^{T}\mathbf{1}=0\\\left\Vert x\right\Vert =1\\y^{T}x=0}}x^{T}Lx$
\For{$i=1,\ldots,n$}
\State $u=\frac{x_{i}}{\sqrt{x_{i}^{2}+y_{i}^{2}}}x+\frac{y_{i}}{\sqrt{x_{i}^{2}+y_{i}^{2}}}y$ 
\State Set $R$ with the $n/2$ vertices with largest $u_j$ and $S$ with the remaining vertices. 
\If{$E(R,S)<E(A,B)$} 
\State $A=R$  
\State $B=S$  
\EndIf 
\EndFor \\
\Return $\{A,B\}$ \\
\end{algorithmic}

\caption{Spectral Bisection Refinement.}
\label{algo:algorithmImprove}
\end{algorithm}

\section{\label{sec:Experimentalresults}Experimental results}

We compared the quality of partitions returned by SB and Algorithm
\ref{algo:algorithm1} on a wide range of graph matrices. The matrices
represents graphs arising in different application domains found in
Matrix Market. Table \ref{tab:tab1} describes the characteristics
of these matrices and the comparison between cut sizes of both algorithms.

\begin{table}[h]
\begin{tabular}{|c|c|c|c|c|c|}
\hline 
\multirow{2}{*}{Matrix} & \multirow{2}{*}{Description} & \multirow{2}{*}{Order} & SB & \multicolumn{2}{c|}{Algorithm \ref{algo:algorithm1}}\tabularnewline
\cline{4-6} 
 &  &  & Cut & Cut & Improv\tabularnewline
\hline 
\hline 
cegb3306 & Structural engineering & 3306 & 18281 & 2421 & \textbf{86\%}\tabularnewline
\hline 
cegb3024 & Structural engineering & 3024 & 19660 & 19534 & 0.6\%\tabularnewline
\hline 
dwt\_1242 & Structural engineering & 1242 & 101 & 72 & \textbf{28\%}\tabularnewline
\hline 
dwt\_2680 & Structural engineering & 2680 & 85 & 85 & 0\%\tabularnewline
\hline 
dwt\_918 & Structural engineering & 918 & 71 & 61 & 14\%\tabularnewline
\hline 
eris1176 & Electrical network & 1176 & 313 & 202 & \textbf{35\%}\tabularnewline
\hline 
bcspwr10 & Power network & 5300 & 44 & 31 & \textbf{29\%}\tabularnewline
\hline 
jagmesh1 & Finite element model & 936s & 50 & 50 & 0\%\tabularnewline
\hline 
jagmesh7 & Finite element model & 1138 & 29 & 28  & 3.4\%\tabularnewline
\hline 
lock2232 & Structural engineering & 2232 & 1008 & 977 & 3\%\tabularnewline
\hline 
lshp1270 & Finite element model  & 1270 & 73 & 73 & 0\%\tabularnewline
\hline 
lshp1882 & Finite element model & 1882 & 89 & 89 & 0\%\tabularnewline
\hline 
commanche\_dual & Structural engineering & 7920 & 46 & 42 & 8.6\%\tabularnewline
\hline 
lshp2614 & Finite element model & 2614 & 105 & 105 & 0\%\tabularnewline
\hline 
lshp3466 & Finite element model & 3466 & 121 & 121 & 0\%\tabularnewline
\hline 
man\_5976 & Structural engineering & 5976 & 55682 & 55391 & 0.5\%\tabularnewline
\hline 
\end{tabular}

\caption{Comparative analysis between SB and Algorithm \ref{algo:algorithm1}.}
\label{tab:tab1}
\end{table}

The last column of Table \ref{tab:tab1} indicates percentage of improvement
of Algorithm \ref{algo:algorithm1} over SB. We highlight the best
results

Next, we compared the quality of partitions for several random graphs
by computing the average gain of Algorithm \ref{algo:algorithm1}
over SB. Here, a random graphs with $n$ vertices follows the Erd\H{o}s\textendash Rényi
model, where an edge is present between two vertices uniformly with
probability $p$. For different combinations of probabilities and
number of vertices, we sampled 1000 random graphs and calculated the
average gain. The experiments discarded graphs that are disconnected.
Table \ref{tab:randomG} shows the resulting ratio of improvement,
where each column corresponds to a given number of vertices $n$ and
each row to a given probability $p$. 

\begin{table}[h]
\begin{tabular}{|c|c|c|c|}
\hline 
$p\backslash n$ & 100 & 500 & 1000\tabularnewline
\hline 
\hline 
0.1 & \textbf{7.68\% } & \textbf{1.87\% } & \textbf{1.01\%} \tabularnewline
\hline 
0.2 & \textbf{4.67\%}  & 0.90\%  & 0.46\% \tabularnewline
\hline 
0.3 & \textbf{3.29\%}  & 0.61\%  & 0.30\% \tabularnewline
\hline 
0.4 & 2.73\%  & 0.64\%  & 0.32\% \tabularnewline
\hline 
0.5 & 2.62\%  & \textbf{1.10\% } & \textbf{0.82\%} \tabularnewline
\hline 
0.6 & 2.91\%  & \textbf{1.20\% } & \textbf{0.68\%} \tabularnewline
\hline 
0.7 & 1.98\%  & 0.39\%  & 0.19\% \tabularnewline
\hline 
0.8 & 1.00\%  & 0.19\%  & 0.09\% \tabularnewline
\hline 
0.9 & 0.57\%  & 0.15\%  & 0.07\% \tabularnewline
\hline 
\end{tabular}

\caption{Average gain for 1000 random graphs}
\label{tab:randomG}
\end{table}

The expected number of edges of these random graphs is $pn(n-1)/2$.
Thus, Table \ref{tab:randomG} indicates that Algorithm \ref{algo:algorithm1}
performs better for sparse graphs than for dense graphs. We highlighted
three best results for each column of Table \ref{tab:randomG}. We
notice that in multilevel algorithms, the coarsest graph is usually
small, with 100 vertices or less. Putting $p=0.1$ we obtain on the
average 495 edges for random graphs with 100 vertices. Table \ref{tab:randomG}
indicates a good improvement ratio for those graphs, with average
of 7.6\%. That suggests that very often the new algorithm provides
better cuts for the initial partition in multilevel algorithms.

\end{document}